\documentclass[12pt]{article}
\usepackage{amsfonts}
\begin{document}
\def\R{\mathbb{R}}
\def\C{\mathbb{C}}
\def\Z{\mathbb{Z}}
\def\N{\mathbb{N}}
\def\Q{\mathbb{Q}}
\def\D{\mathbb{D}}
\def\T{\mathbb{T}}
\def\hb{\hfil \break}
\def\ni{\noindent}
\def\i{\indent}
\def\a{\alpha}
\def\b{\beta}
\def\e{\epsilon}
\def\d{\delta}
\def\De{\Delta}
\def\g{\gamma}
\def\qq{\qquad}
\def\L{\Lambda}
\def\E{\cal E}
\def\G{\Gamma}
\def\F{\cal F}
\def\K{\cal K}
\def\A{\cal A}
\def\B{\cal B}
\def\M{\cal M}
\def\P{\cal P}
\def\Om{\Omega}
\def\om{\omega}
\def\s{\sigma}
\def\t{\theta}
\def\th{\theta}
\def\Th{\Theta}
\def\z{\zeta}
\def\p{\phi}
\def\P{\Phi}
\def\m{\mu}
\def\n{\nu}
\def\l{\lambda}
\def\Si{\Sigma}
\def\q{\quad}
\def\qq{\qquad}
\def\half{\frac{1}{2}}
\def\hb{\hfil \break}
\def\half{\frac{1}{2}}
\def\pa{\partial}
\def\r{\rho}
\begin{center}
{\bf SURVEY ARTICLE}
\end{center}
\begin{center}
{\bf RIESZ MEANS AND BEURLING MOVING AVERAGES}
\end{center}
\begin{center}
{\bf N. H. BINGHAM}
\end{center}

\ni {\bf Abstract}.  We survey the interplay between the Riesz means and Beurling moving averages of the title, obtaining Abelian and Tauberian results relating different Riesz means (or Beurling moving averages) whose defining functions have comparable growth.  The motivation includes strong limit theorems in probability theory. \\
\ni {\it Keywords}: Riesz means, Beurling's Tauberian theorem, moving average, regular variation, laws of large numbers, Dirichlet series, Prime Number Theorem \\
\ni {\it MSC Classification numbers}: 40E05, 60F15, 26A15. \\

\ni {\bf 1.  Introduction} \\
\i Riesz (or typical) means have proved very useful since their introduction by Marcel Riesz in 1909.  He used them for summing Dirichlet series, and so providing analytic continuation; see e.g. the Cambridge Tract by Hardy and Riesz [HarR], and the monograph by Chandrasekharan and Minakshisundaram [ChaM III].  A Riesz mean $R(\l, k)$ has an {\it order} $k > 0$ and a {\it type} $\l$, where $\l: {\R}_+ \to {\R}_+$ is a function increasing to infinity; we write ${\l}_n$ for $\l(n)$, and regard $\l$ as specified by these values.  We confine attention here to order $k = 1$, and abbreviate $R(\l, 1)$ to $R(\l)$. \\
\i Such Riesz means were used extensively by Karamata [Kar2] in the 1930s.  Writing
$$
s \to c \q R(\l) \qq \hbox{for} \qq \frac{1}{\l(x)} \int_0^x s(u) d \l(u) \to c \q (x \to \infty), \eqno(R(\l))
$$
he gave necessary and sufficient Tauberian conditions for such Riesz convergence to imply ordinary convergence $s(x) \to c$ as $x \to \infty$. \\
\i The whole field of Tauberian theory was revolutionised by Wiener in 1932 [Wie1-3], who passed from special kernels to general kernels.  The Wiener Tauberian theorem deals with convolutions (of a function with a kernel), which may be written multiplicatively (using the group $({\R}_+, \times)$, and Mellin transforms) or additively (using the group $(\R, +)$ and Fourier transforms); for textbook accounts, see [Wie2], Hardy [Har4 XII], Korevaar [Kor, II], Widder [Wid, V]. \\
\i Not all important Tauberian questions lend themselves to treatment in such a convolution framework, most notably that of the {\it Borel} summability method and its relatives ([Har4 VIII, IX], [Kor VI]).  In unpublished lectures of 1957, Beurling gave an extension of the Wiener Tauberian theorem, allowing one to pass from statements of the form
$$
\int_{\R} S(u) K \Bigl( \frac{x - u}{\p(x)} \Bigr) \frac{du}{\p(x)} \to c \int_{\R} K(u) du \qq (x \to \infty)
$$
to similar statements with the kernel $K$ replaced by $H$; here $H, K$ are in $L_1(\R)$, $K$ is a {\it Wiener kernel} (its Fourier transform has no real zeros), $S$ is bounded, and the function $\p(x): {\R}_+ \to {\R}_+$ is $o(x)$ as $x \to \infty$, measurable, and satisfies
$$
\p(x + t \p(x))/\p(x) \to 1 \qq (x \to \infty) \qq \forall t \in \R; \eqno(BSV)
$$
call such functions $\p$ {\it Beurling slowly varying} (notation: $\p \in BSV$), and call the left-hand side here the {\it Beurling convolution} with respect to $\p \in BSV$.  Then Wiener's Tauberian theorem extends to this new setting, to give {\it Beurling's Tauberian theorem}, which reduces to it in the special case $\p \equiv 1$; for a textbook account, see [Kor IV.11]. \\
\i If one takes for $H$ the indicator function of an interval, the conclusion of Beurling's Tauberian theorem takes the form
$$
\frac{1}{c \p(x)} \int_x^{x + c \p(x)} S(u) du \to s \qq (x \to \infty) \qq \forall c > 0. \eqno(BMA_{\p})
$$
We call the left-hand side here the {\it Beurling moving average} w.r.t. $\p$.  We shall be principally interested in sequences $(s_n)$; writing $S(x) := \sum_{n \leq x} s_n$, the above becomes
$$
\frac{1}{c \p(x)} \sum_{x < n \leq x + c \p(x)} s_n \to s \qq (x \to \infty) \qq \forall c > 0. \eqno(BMA_{\p})
$$
\i Such Beurling moving averages have many uses -- for example, in probability theory in connection with strong laws of large numbers -- and have been studied from that point of view in e.g. [Bin1-7], [BinGo1,2], [BinGa], [BinM], [BinR], [BinT].  They have also been studied recently from the point of view of Beurling slow and regular variation (see [BinO1-4]).  This is related to the important subject of Karamata regular variation [Kar1], to which we return below; for textbook accounts see e.g. [BinGoT] (BGT for brevity), [Kor IV]. \\
\i Our aim here is to study the inter-relationship between different Beurling moving averages, as $\p$ varies (increases or decreases).  As $\p$ {\it increases} in $(BMA_{\p})$, more averaging takes place, and this suggests that matters improve and convergence is preserved.  It has been known since [Bin1] that Beurling averages are special cases of Riesz means.  For Riesz means $R(\l)$, matters improve as $\l$ {\it decreases}.  Furthermore, $R(\l)$ and $R(\m)$ are equivalent if $\log \l$ and $\log \m$ are {\it comparable} ($0 < \liminf \log \l /\log \m \leq\limsup \log \l /\log \m < \infty$).  This is the second consistency theorem for Riesz means ([HarR Th. 17] for the case $\l \mapsto \log \l$, [Har2, \S 8], [Har3] in general; [ChaM II]).  This becomes transparent in our context of Beurling moving averages, as we shall see in \S 2 below.   \\

\ni {\bf 2.  Results} \\
\i We shall need {\it local uniformity} -- that is, that $(BSV)$ holds uniformly on compact $t$-sets.  For this, we need to impose a condition on $\p$, and we have several choices here. \\
1.  We can take $\p$ {\it continuous}.  Then uniformity holds by Bloom's theorem of 1976 ([Blo]; BGT, \S 2.11; [BinO1]). \\
2.  We can take $\p$ monotone [BinO1]. \\
3.  We can take $\p$ measurable and with the Darboux property (intermediate-value property): [BinO1]. \\
4.  We can take $\p$ with the Baire and Darboux properties [BinO1]. \\
These are all quite weak; indeed, continuity and monotonicity usually hold from the context in practice. \\
\i When one has local uniformity in $(BSV)$, one writes $\p \in SN$, and calls $\p$ {\it self-neglecting}; we restrict attention to such $\p$ below. \\
\i We have mentioned the relevance of Karamata's theory of regular variation in \S 1 above (see e.g. BGT Ch. 1).  Here the main result is the {\it uniform convergence theorem} (UCT): if
$$
f(\l x)/f(x) \to g(\l) \qq (x \to \infty) \qq \forall \l > 0, \eqno(K)
$$
then the convergence is uniform on compact $\l$-sets in $(0,\infty)$ if $f$ has the Baire property (`is Baire') or is measurable, and the limit $g$ has the form $g(\l) \equiv {\l}^{\rho}$ for some $\rho$ (the {\it characterisation theorem}).  There is a similar UCT in the Bojanic-Karamata/de Haan theory of regular variation (see e.g. BGT Ch. 3, and for the Bojanic-Karamata contribution, [BojK]).  Whereas Beurling slow variation originally appeared as a minor topic in Karamata theory (see e.g. BGT \S 2.11), it is now known that Beurling slow and regular variation {\it include} the Karamata and Bojanic-Karamata/de Haan theories [BinO1].  It thus shares their feature, of having a {\it representation theorem}: one can give an explicit description of all functions in $(BMA_{\p})$, in $SN$, etc. Our results follow simply from this.  The functions in $SN$ have representations of the form \\
$$
\phi(x) = c(x) \int_0^x \e(u) du
$$
with $c(.) \to c \in (0,\infty)$, $\e(.) \to 0$ and the product positive (BGT \S 2.11; [BinO1]); by uniformity in $(BMA_{\p})$, one may take $\p(x) = c \int_0^x \e(u) du$, so ${\p}'(x) = c \e(x) \to 0$. \\
\i We quote the following results from [BinGo2, Th. 1,2] and [BinGo1, Th. 3] respectively.  The first links the two themes of our title, the second is a representation theorem. \\

\ni {\bf Theorem BG1}.  For $\p \in SN$, write
$$
\P(x) := \int_1^x du/\p(u).
$$
Then $(BMA_{\p})$ for a sequence $(s_n)$ implies $(R(\l))$ and $(R(\m))$, where
$$
\l(x) := \p(x) \exp \{ \P(x) \}, \qq \m(x) := \exp \{ \P(x) \}.
$$

\ni {\it Note}.  In the second consistency theorem, what is crucial is $\log \l$.  This is $\log \p + \P$, and as $\p(x) = o(x)$, $\log \p$ is $o(\log x)$.  But $\P(x) = \int_1^x du/\p(u)$ is of order greater than $\int_1^x du/u = \log x$, so $\log \l \sim \log \mu = \P$.  Thus comparability of two $\p$-functions gives comparability of their $\P$-functions on integrating, so of their ($\log \mu$ or) $\log \l$ functions, so {\it equivalence} of the Riesz means.  Further, convergence improves as $\p$ increases, that is, $\l$ decreases, as expected in \S 1. \\

\ni {\bf Theorem BG2}.  The function $U$ satisfies $(BMA_{\p})$ iff it is of the form
$$
U(x) = const + o(\p(x)) + c \int_1^x (1 + \eta(u))du \q \hbox{with} \q \eta(.) \to 0.
$$

\i This provides a transparent proof of the following result, which links Beurling moving averages $BMA_{\p}$ with different self-neglecting functions $\p$. \\

\ni {\bf Theorem}.  If $\p, \psi \in SN$ with $\p \leq \psi$: \\
(i) $U \to c$ $(BMA_{\p})$ implies $U \to c$ $(BMA_{\psi})$. \\
(ii) The converse holds iff there exists a function $\d \to 0$ such that $V := U + \d$ satisfies $(BMA_{\p})$. \\

\ni {\it Proof}.  The direct assertion (i) just says that $o(\p)$ is $o(\psi)$; this is Abelian.  The converse assertion (ii) says that {\it if} there is an $o(1)$-perturbation $V = U + \d$ of $U$ satisfying $(BMA)_{\p}$, then so does $U$.  This is just linearity in $(BMA)_{\p}$.  Also, this perturbation condition is necessary as well as sufficient (if the conclusion holds, the perturbation condition holds with $\d = 0$), and is thus a {\it necessary and sufficient Tauberian condition}. // \\

\i This result, though nearly immediate given the representation theorem Th. BG2, has distinguished antecedents.  For the Abelian result (i), the background is (what became known as) the second consistency theorem for Riesz means (\S 1); this emerged in work of Riesz in 1909 and 1916 [Rie1,2], Hardy [Har2,3] in 1910 and 1915, [HarR Th. 17] in 1915.  For the Tauberian result (ii), the origin is Hardy's result of 1904 ([Har1, \S 6, p. 40]; [Har2, 55; Comments, 83]; [Har4, Th. 149]), that Ces\`aro convergence {\it with rate} $o(1/\sqrt{n})$ gives Borel convergence.  (This is extended to general order $k$ for Borel's integral method by Hardy and Littlewood in 1916 [HarL, \S 3].)  It was further extended to Euler convergence by Knopp in 1923 [Kno]. \\
\i Hardy and Littlewood introduced a family of summability methods in [HarL], for the purpose of analytic continuation by power series [Har4, VIII].  They called these the $(e,c)$ methods, though the term {\it Valiron methods} is now used; there is a textbook treatment in [Har4, IX, 9.10, 9.11] in 1949.  Also in 1949, Meyer-K\"onig [MeyK] gave an extended treatment of the circle methods ({\it Kreisverfahren}: methods of Euler, Borel, Valiron, Taylor and Meyer-K\"onig; there is a modern textbook account in [Kor, VI]); unfortunately, this came too late for inclusion in [Har4].  The Riesz method linked to the Borel, Euler and circle methods is $R(e^{\sqrt{n}})$.  This is equivalent, both to the Beurling moving average $BM_{\sqrt{n}}$ and to Ces\`aro convergence with rate $o(1/\sqrt{n})$ of some $o(1)$-perturbation as in (ii) of the Theorem [Bin2,3].  The Abelian result that convergence under $R(e^{\sqrt{n}})$ implies Euler (and so Borel) and Valiron convergence follows from Knopp's result in the Euler case and work of Hyslop [Hys] in the Valiron case, Parameswaran [Par] in the circle case; see also [Bin2,3].  The Tauberian converse is proved for $s_n$ non-negative (equivalently, bounded below: $s_n = O_L(1)$) by Tenenbaum [Ten].  It is given under a best-possible one-sided Tauberian condition in [BinGo1]; cf. [Bin1], [BinT].  That all these methods are equivalent for {\it bounded} sequences is due to Meyer-K\"onig [MeyK]; for a simple proof in the Euler-Borel case, see Jurkat [Jur].  Equivalence of all these methods for bounded sequences follows from Beurling's Tauberian theorem (\S 1; see [Kor, IV.11]), developed by Beurling for this purpose. \\
\i For $0 < \b < 1$, equivalence of the Riesz mean $R(\exp \{ n^{1 - \b} \})$ with $BMA_{n^{\b}}$ and with perturbed Ces\`aro convergence with rate $o(1/n^{1-\b})$ is in [BinT, Th. 3].  One can extend this to include $\b = 1$, by using $R(\exp \{ \int_1^n dx/x^{\b})$) (see \S 3.4 below).  Part of this result was obtained much earlier by Wang [Wan]. \\
\i The one-sided Tauberian condition above derives from work of Erd\H{o}s and Karamata [ErdK] in 1956 and Bojanic and Karamata [BojK] in 1963; see BGT \S 3.8 for a textbook account. \\

\ni {\bf 3.  Remarks}. \\
 1. {\it Riesz means and their summability}.  \\
 \i The theory of Riesz means owes its development (by Riesz, 1909-1916: [Rie1-3]) to the intense interest in Dirichlet series in general and the Riemann zeta function in particular in the early 1900s, following the proof of the Prime Number Theorem in 1896.  For background, see e.g. Landau [Lan, Sechstes Buch], [HarR], [BohC, I.5,6].  For example, the Riemann zeta function cannot be extended by any Ces\`aro method of summability [Har2, \S 6].  However, the {\it alternating} zeta function $\sum_1^{\infty} (-)^{n-1}/n^s$ is summable by Ces\`aro means of order $k$ if $\sigma := Re \ s > -k$ (Bohr [Boh1, II.3 Th. I]; [HarR, 20, 24]), and as $\sum_1^{\infty} (-)^{n-1}/n^s = \zeta(s) (1 - 2^{1-s})$, this can be used to continue $\zeta$.  See [HarR IV] for general results, and [HarR VI] for Abelian and Tauberian theorems.    \\
2.  {\it The logarithmic method}.  \\
\i One writes
$$
\frac{1}{\log n} \sum_0^n s_k/(k+1) \to s \qq (n \to \infty)
$$
as
$$
s_n \to s \qq (\ell),
$$
and says that $s_n$ converges to $s$ under the {\it logarithmic method} of summability.  Then $\ell$ is the Riesz method $R(\log n)$ (see e.g. [Har4, Th. 37]).  As $\log n$ grows more slowly than $n$, this method is {\it less} restrictive than the Ces\`aro method $C = C_1 = R(n)$, while $R(e^{\sqrt{n}})$ is {\it more} restrictive.  To summarise:
$$
R(\exp \{ \sqrt{n} \}) \subset R(n) \subset R(\log n).
$$
The logarithmic method has recently been studied in detail by the author and Gashi [BinGa].  One obtains moving averages there too, but these are {\it not} of the Beurling type of \S \S 1,2, [BinGo1,2] etc.    \\
3.  {\it Regular variation}.  \\
\i As noted above, Beurling regular variation [BinO4] contains both the original Karamata regular variation and its later development by Bojanic and Karamata and by de Haan.  A structural feature of {\it any} theory of regular variation is a {\it representation theorem}.  As the simple proof of our Theorem clearly shows, it is this that is crucial here.  It is precisely the lack of any awareness of this aspect -- that once one can identify regular variation in some form, a complete and explicit description of the function in question is available -- that delayed the solution of the problems considered above.  \\
4.  {\it Laws of Large Numbers (LLNs)}.  \\
\i The central result of probability theory is Kolmogorov's strong law of large numbers (SLLN) of 1933 [Kol]: for $X, X_1, \ldots, X_n, \ldots$ independent and identically distributed (iid) random variables, the sample mean $S_n/n := \frac{1}{n} \sum_1^n X_k$ converges with probability one (almost surely, a.s.) to some constant $c$ iff $X \in L_1$, i.e. $E[|X|] < \infty$; then $c = E[X]$, the population mean.  The summability method here is the Ces\`aro method $C = C_1$; the moment condition here is membership of $L_1$. \\
\i This result has variants of many kinds; the one that concerns us here is alternatives in which both the summability method and the integrability condition are changed.  Both are {\it strengthened} if one works with
$$
R_p := R(\exp \{\int_1^n dx/x^{1/p} \}) \q (p \geq 1);
$$
it turns out that the corresponding integrability condition is membership of $L_p$ (see \S 3.7 below).  The most important case is $p = 2$ discussed above, where the method is $R_2$, equivalently $R(e^{\sqrt{n}})$, and the integrability condition is existence of a {\it variance} (as distinct from Kolmogorov's SLLN for $p = 1$, which is about {\it means}).  It turns out that the Euler and Borel (and Valiron) methods are interchangeable with $R_2$ here; for details and references, see e.g. [BinT].  But one is not restricted to {\it powers} here, as is shown in [BinGo2, \S 3].  There, the summability method is $R(\l)$, with $\l(x) := \exp \{ \int_1^x du/\p(u) \}$ with $\p \in SN$; the integrability condition is$E[\psi(|X|)] < \infty$ with $\psi := {\p}^{\leftarrow}$ and ${\psi}'$ assumed of {\it bounded increase} (roughly, polynomial growth).  This restriction is necessary for probabilistic reasons (a.s. invariance v. non-invariance, Erd\H{os}-R\'enyi law, etc.).  The proof hinges on {\it maximal inequalities}, as usual with strong limit theorems in probability theory.  Again, one sees that, with two comparable functions $\p$ (or $\psi := {\p}^{\leftarrow}$), it is the bigger $\psi$ (smaller $\p$, bigger $\l := \exp \{ \int_1^. 1/\p \}$) that is the more restrictive. \\
\i By contrast, both summability and integrability are {\it weakened} when we pass from $C$ to the logarithmic method.  It turns out that this case is very interesting and important, because of its connection with the almost sure central limit theorem (ASCLT).  We refer to [BinR], [Bin6], [BinGa] for details and references here.  The logarithmic method is also important in analytic number theory, in connection with the Prime Number Theorem and with densities of sets of integers; again, see [BinGa] for details.   \\
5.  {\it CLT; LIL, LSL}.  \\
\i The second main pillar of probability theory is the central limit theorem (CLT: the folklore names of LLN and CLT are the Law of Averages and the Law of Errors).  This too has a counterpart in this context, due to Embrechts and Maejima [EmbM].  Intermediate between the two is the law of the iterated logarithm (LIL).  This applies to the Ces\`aro method, $C = R_1$; its counterpart for $R_2$ contains only one logarithm, and is called the law of the single logarithm (LSL); see e.g. [BinM], [Bin5] for details and references.  For $R_p$, Gut, Jonsson and Stadtm\"uller [GutJS] obtain a result which is, in a sense, {\it between LIL and LSL}. \\
6.  {\it Tauberian theorems between Riesz means}.  \\
\i We mention here the work of Jurkat, Kratz and Peyerimhoff [JurKP].  Here, given three Riesz means $R({\l}_i, k_i)$ ($i = 1,2,3$), if two Riesz means of a function have given bounds, best-possible bounds are given for the third.  These results are used in [Bin5] to show that for $p > 1$, the Riesz mean $R_p$ is strictly more restrictive than the Ces\`aro mean $C_{1/p}$. \\
7.  {\it Discontinuity across} $p = 1$. \\
\i The Ces\`aro methods $C_{\a}$ for $\a \geq 1$ are {\it equivalent} in the context above of LLNs (to each other and to the Abel method $A$):
$$
X \in L_1 \q \Leftrightarrow \q X_n \to \mu \q a.s. \q (C_1) \q \Leftrightarrow X_n \to \mu \q a.s. \q (C_{\a}) \q (\a \geq 1).
$$
By contrast, for $\a \leq 1$ ($\a = 1/p$ above) the integrability required is the more stringent $L_{1/\a}$:
$$
X \in L_{1/\a} \q \Leftrightarrow \q X_n \to \mu \q a.s. \q (R_{1/\a}) \q \Leftrightarrow X_n \to \mu \q a.s. \q (C_{\a}) \q (0 < \a \leq 1);
$$
see [Bin5].  This striking discontinuity reflects the distinguished role played by the Kolmogorov SLLN among results of this type.  This is seen also in the different nature of the LLNs for the logarithmic method [BinGa], where one assumes {\it less} integrability than in the SLLN, rather than more as above.  \\
\i Despite this, there are two results which handle $L_p$ with $p \leq 1$ and $p \geq 1$ together (though there is a difference even there: one centres at the mean where the mean exists): the Marcinkiewicz-Zygmund law and the Baum-Katz law [MarZ], [BauK].  For a textbook account, see [Gut, 6.7, 6.12.1]; for a survey of extensions to the strong law generally, see [Bin4]. \\
\ni 8.  {\it Self-neglecting and self-equivarying functions}. \\
\i The self-neglecting functions $\p \in SN$ of \S 1 have $\p(x) = o(x)$ (by definition, or from the representation).  It turns out that one may usefully extend this to $\p(x) = O(x)$, but no further; the resulting functions are called {\it self-equivarying} functions, forming the class $SE$.  It is this extension from $SN$ to $SE$ that enables the Karamata and Bojanic-Karamata/de Haan theories of regular variation to be subsumed within that of Beurling regular variation [BinO1], [Ost1,2].  But this shows why the logarithmic method lies outside the Beurling framework: $\l(x) = \log x = \exp \{ \int_1^x du/\p(u) \}$ would give $\p(x) = x \log x$, which is not $O(x)$. \\
9.  {\it Systems of kernels}. \\
\i We saw in the Beurling Tauberian theorem the use of $H := I_{[-c,0]}$.  This cannot be used as a Wiener kernel, as its Fourier transform has zeros $2n \pi/c$.  But for two such $H_i$ with $c_i$ incommensurable, there are no {\it common} zeros, and in the Wiener Tauberian theory (and its Beurling extension) this suffices.  See [BinI2] for a case in point (the motivation is analytic number theory, in [BinI1]).  For by Kronecker's theorem (see e.g. [HarW, XXIII, Th. 438]) the additive group generated by $c_1, c_2$ is then dense in $\R$.  This bears on {\it quantifier weakening}: the quantifier $\forall$ in $BMA_{\p}$ may be weakened to just {\it two} such incommensurable $c$s.  This is a common theme in regular variation; see BGT Th. 3.2.5, and the recent [BinO2]. \\

\ni {\bf Postscript}.  This paper was written in the summer of 2014, shortly before the centenary of the outbreak of World War I.  It will probably appear in 2015, the centenary of the Hardy-Riesz Tract (from the Preface, by Hardy, 19 May 1915: `The publication of this tract has been delayed by a variety of causes, and I am now compelled to issue it without Dr Riesz's help in the final correction of the proofs': Riesz, who was Hungarian, was in Stockholm).  We recall their unforgettable dedication: \\
\i Mathematicis quotquot ubique sunt operum societatem nunc diremptam mox ut optare licet redintegraturis d.d.d. auctores hostes idemque amici.\\
\i [To mathematicians, however many and wherever they may be: that they may soon again take up, as is to be hoped, the fellowship of their work now disrupted, we the authors, enemies and yet also friends, present and dedicate (this book).] \\

\ni {\bf Acknowledgements}.  It is a pleasure to thank three valued collaborators, one old and two more recent.  The old one is Charles Goldie, with whom I worked on such things in the 1980s.  The more recent ones are Bujar Gashi and Adam Ostaszewski, who rekindled my interest in the Riesz and Beurling aspects respectively. \\

\begin{center}
{\bf References}
\end{center}
\ni [BauK] L. E. Baum and M. Katz, Convergence rates in the law of large numbers.  {\sl Trans. Amer. Math. Soc.} {\bf 120} (1965), 108-123. \\
\ni [Bin1] N. H. Bingham, Tauberian theorems and the central limit theorem.  {\sl Ann. Probab.} {\bf 9} (1981), 221-231. \\
\ni [Bin2] N. H. Bingham, On Euler and Borel summability.  {\sl J. London Math. Soc.} (2) {\bf 29} (1984), 141-146. \\
\ni [Bin3] N. H. Bingham, On Valiron and circle convergence.  {\sl Math. Z.} {\bf 186} (1984), 273-286.  \\
\ni [Bin4] N. H. Bingham,  Extensions of the strong law.  {\sl Analytic and Geometric Stochastics} (ed. D.G. Kendall), 27-36.  Supplement, Adv. Appl. Probability
(G.E.H. Reuter Festschrift), 1986. \\
\ni [Bin5] N. H. Bingham, Moving averages.  {\sl Almost Everywhere Convergence I}  (ed. G.A. Edgar \& L. Sucheston) 131-144, Academic Press, 1989. \\
\ni [Bin6] N. H. Bingham,  The strong arc-sine law in higher dimensions.  {\sl Convergence in Ergodic Theory and Probability} (ed. V. Bergelson, P.
March \& J. M. Rosenblatt) 111-116, Walter de Gruyter, Berlin - New York, 1996. \\
\ni [Bin7] N. H. Bingham, Hardy, Littlewood and probability.  {\sl Bull. London Math. Soc.}, to appear. \\
\ni [BinGa] N. H. Bingham and Bujar Gashi, Logarithmic moving averages.  {\sl J. Math. Anal. Appl.} {\bf 421} (2015), 1790-1802. \\
\ni [BinGo1] N. H. Bingham and C.M. Goldie,  On one-sided Tauberian conditions. {\sl Analysis} {\bf 3} (1983), 159-188. \\
\ni [BinGo2] N. H. Bingham and C.M. Goldie,  Riesz means and self-neglecting functions. {\sl Math. Z.} {\bf 199} (1988), 443-454. \\
\ni [BinGoT] N. H. Bingham, C. M. Goldie and J. L. Teugels, {\sl Regular variation}.  Encycl. Math. Appl. {\bf 27}, Cambridge Univ. Press, 1987 (2nd ed. 1989). \\
\ni [BinI1] N. H. Bingham and A. Inoue, Abelian, Tauberian and Mercerian theorems for arithmetic sums.  {\sl J. Math. Anal. Appl.} {\bf 250} (2000), 465-493. \\
\ni [BinI2] N. H. Bingham and A. Inoue, Tauberian and Mercerian theorems for systems of kernels.  {\sl J. Math. Anal. Appl.} {\bf 252} (2000), 177-197.  \\
\ni [BinM] N. H. Bingham and M. Maejima,  Summability methods and almost-sure convergence. {\sl Z. Wahrschein.} {\bf 68} (1985), 383-392  \\
\ni [BinO1] N. H. Bingham and A. J. Ostaszewski, Beurling slow and regular variation.  {\sl Transactions of the London Mathematical Society} {\bf 1} (2014), 29-56. \\
%(see also Uniformity and self-neglecting functions, arXiv:1301.5894; II. Beurling regular variation and the class $\Gamma$, arXiv:1307.5305). \\
\ni [BinO2] N. H. Bingham and A. J. Ostaszewski, Additivity, subadditivity and linearity: automatic continuity and quantifier weakening.  arXiv:1405.3948; submitted. \\
\ni [BinO3] N. H. Bingham and A. J. Ostaszewski, Cauchy's functional equation and extensions: Goldie's equation and inequality, the Golab-Schinzel equation and Beurling's equation.  arXiv:1405.3947; submitted. \\
\ni [BinO4] N. H. Bingham and A. J. Ostaszewski, Beurling moving averages and approximate homomorphisms.  arXiv:1407.4093; submitted. \\
\ni [BinR] N. H. Bingham and L.C.G. Rogers,  Summability methods and almost-sure convergence. {\sl Almost Everywhere Convergence II}  (ed. A. Bellow \& R.L. Jones), Academic Press (1991) 69-83. \\
\ni [BinT] N. H. Bingham and G. Tenenbaum,  Riesz and Valiron means and fractional moments. {\sl Math. Proc. Cambridge Phil. Soc.} {\bf 99} (1986), 143-149. \\
\ni [Blo] S. Bloom, A characterization of B-slowly varying functions.  {\sl Proc. Amer. Math. Soc.} {\bf 54} (1976), 243-250. \\
\ni [Boh1] H. Bohr, Contributions to the theory of Dirichlet series.  Paper S1 in [Boh2] (translation of Bohr's thesis of 1910, in Danish). \\
\ni [Boh2] H. Bohr, {\sl Collected mathematical works}, Vol. III, Danish Math. Soc., Copenhagen, 1952. \\
\ni [BohC] H. Bohr and H. Cram\'er, Die neuere Entwicklung der analytischen Zahlentheorie.  {\sl Encycl. math. Wiss.} II, C 8 (1923), 722-849 (reprinted as H in [Boh2] and Harald Cram\'er,
{\sl Collected Works} Vol. I, Springer, 1994, 289-416). \\
\ni [BojK] R. Bojanic and J. Karamata, On a class of functions of regular asymptotic behaviour.  {\sl Math. Research Center Tech. Reports} {\bf 436} (1963), Univ. Wisconsin, Madison (reprinted, [Kar3, 545-569]). \\
% {\sl Jovan Karamata: Selected papers}, Zavod za Ud\v{z}benike, Beograd, 2009, 545-569). \\
\ni [ChaM] K. Chandrasekharan and S. Minakshisundaram, {\sl Typical means}, Oxford Univ. Press, 1952.\\
\ni [EmbM] P. Embrechts and M. Maejima, The central limit theorem for summability methods of i.i.d. random variables.  {\sl Z. Wahrschein.} {\bf 68} (1984), 191-204.  \\
\ni [ErdK] P. Erd\H{o}s and J. Karamata, Sur la majorabilit\'e $C$ des suites de nombres r\'eels.  {\sl Publ. Inst. Math. (Beograd)} {\bf 10} (1956), 37-52.   \\
\ni [Gut] A. Gut, {\sl Probability: a graduate course}.  Springer, 2005. \\
\ni [GutJS] A. Gut, F. Jonsson and U. Stadtm\"uller, Between the LIL and the LSL.  {\sl Bernoulli} {\bf 16} (2010), 1-22.   \\
\ni [Haa] L. de Haan, {\sl On regular variation and its application to the weak convergence of sample extremes}.  Math. Centre Tracts {\bf 32}, Mathematisch Centrum, Amsterdam, 1970. \\
\ni [Har1] G. H. Hardy, Researches in the theory of divergent series and divergent integrals.  {\sl Quart. J. Math.} {\bf 35} (1904), 22-66 (reprinted, [Har5, 85-114]). \\
\ni [Har2] G. H. Hardy, Theorems relating to the summability and convergence of slowly oscillating series.  {\sl Proc. London Math. Soc.} {\bf 8} (1910), 301-320 (reprinted, [Har5], 291-312). \\
\ni [Har3] G. H. Hardy, The second theorem of consistency for summable series.  {\sl Proc. London Math. Soc.} {\bf 15} (1915), 72-88 (reprinted, [Har5, 588-605]). \\
\ni [Har4] G. H. Hardy, {\sl Divergent series}, Oxford Univ. Press, 1949. \\
\ni [Har5] G. H. Hardy, {\sl Collected papers Vol. VI: Theory of series}, Oxford Univ. Press, 1974. \\
\ni [HarL] G. H. Hardy and J. E. Littlewood, Theorems concerning the summability of series by Borel's exponential method.  {\sl Rend. Cir. Mat. Palermo} {\bf 41} (1916), 36-53 (reprinted, [Har5, 609-628]). \\
\ni [HarR] G. H. Hardy and M. Riesz, {\sl The general theory of Dirichlet's series}, Cambridge Tracts in Math. {\bf 19}, Cambridge Univ. Press, Cambridge, 1915 (reprinted 1964, Stechert-Hafner). \\
\ni [HarW] G. H. Hardy and E. M. Wright, {\sl An introduction to the theory of numbers}, 6th ed. (rev. D. R. Heath-Brown and J. H. Silverman), Oxford Univ. Press, 2008. \\
\ni [Hys] J. M. Hyslop, On the summability of series by a method of Valiron.  {\sl Proc. London Math. Soc.} {\bf 41} (1936), 243-256. \\
\ni [Jur] W. B. Jurkat, Ein funktiontheoretischer Beweis f\"ur O-Taubers\"atze bei den Verfahren von Borel und Euler-Knopp.  {\sl Arch. Math. (Basel)} {\bf 7} (1956), 278-283. \\
\ni [JurKP] W. B. Jurkat, W. Kratz and A. Peyerimhoff, The Tauberian theorems which inter-relate different Riesz means.  {\sl J. Approximation Th.} {\bf 13} (1975), 2235-266. \\
\ni [Kar1] J. Karamata, Sur un mode de crossiance r\'eguli\`ere des fonctions.  {\sl Mathematica (Cluj)} {\bf 4} (1930), 38-53 (reprinted, [Kar3, 61-79]). \\
\ni [Kar2] J. Karamata, {\sl Sur les th\'eor\`emes inverses des proc\'ed\'es de sommabilit\'e}, Actualit\'es Scientifiques et Industrielles {\bf 450}, Hermann, Paris, 1937. \\
\ni [Kar3] J. Karamata, {\sl Selected papers} (ed. V. Mari\'c).  Zavod za ud\v{z}benike/Math. Inst. Serbian Acad. Sci., Belgrade, 2009. \\
\ni [Kno] K. Knopp, \"Uber das Eulersche Summierungsverfahren II.  {\sl Math. Z.} {\bf 18} (1923), 125-156. \\
\ni [Kol] A. N. Kolmogorov, {\sl Grundbegriffe der Wahrscheinlichkeitsrechnung}.  Erg. Math. {\bf 3}, Springer, 1933 (translated as {\sl Foundations of the theory of probability}, Chelsea, New york, 1956). \\
\ni [Kor] J. Korevaar, {\sl Tauberian theory: A century of developments}.  Grundl. math. Wiss. {\bf 329}, Springer, 2004. \\
\ni [Lan] E. Landau, {\sl Handbuch der Lehre von der Verteilung der Primzahlen}.  Teubner, Leipzig, 1909 (reprinted, Chelsea, New York, 1953). \\
\ni [MarZ] J. Marcinkiewicz and A. Zygmund, Sur les fonctions ind\'ependantes.  {\sl Fund. Math.} {\bf 29} (1937), 60-90 (reprinted, {\sl J\'ozef Marcinkiewicz, Collected papers} (ed. A. Zygmund), PWN, Warszawa, 1964, p.233-259). \\
\ni [MeyK] W. Meyer-K\"onig, Untersuchungen \"uber einige verwandte Limitierungsverfahren.  {\sl Math. Z.} {\bf 52}, (1949), 257-304. \\
\ni [Ost1] A. J. Ostaszewski, Beurling regular variation, Bloom dichotomy and the Golab-Schinzel functional equation.  {\sl Aequationes Math.} (2014), ISSN 0001-9054.  \\
\ni [Ost2] A. J. Ostaszewski, Homomorphisms from functional equations: the Goldie equation II.  arXiv:1407.4089. \\
\ni [Par] M. R. Parameswaran, On summability functions for the circle family of summability methods.  {\sl Proc. Nat. Inst. Sci. India A} {\bf 25} (1959), 171-175. \\
\ni [Rie1] M. Riesz, Sur la sommation des s\'eries de Dirichlet.  {\sl Comptes Rendus Acad. Sci. Paris} {\bf 149} (1909), 18-21 (reprinted in [Rie3, 59-61]).  \\
\ni [Rie2] M. Riesz, Ein Konvergenzsatz f\"ur Dirichletsche Reihen.  {\sl Acta Math.} {\bf 40} (1916), 349-361 (reprinted in [Rie3, 182-194]). \\
\ni [Rie3] M. Riesz, {\sl Collected papers}, Springer, 1988 (reprinted 2013). \\
\ni [Ten] G. Tenenbaum, Sur le proc\'ed\'e de sommation de Borel et la r\'epartition du nombre des facteurs premiers des entiers.  {\sl Enseignement math.} {\bf 26} (1981), 225-245. \\
\ni [Wan] F. T. Wang, A note on Riesz summability of the type $e^{n^{\alpha}}$.  {\sl Bull. Amer. Math. Soc.} {\bf 50} (1944), 417-419. \\
\ni [Wid] D. V. Widder, {\sl The Laplace transform}, Princeton Univ. Press, 1941. \\
\ni [Wie1] N. Wiener, Tauberian theorems, {\sl Ann. Math.} {\bf 33} (1932), 1-100 (reprinted, {\sl Generalized harmonic analysis and Tauberian theorems}, MIT Press, Cambridge MA, 1964 and [Wie3, 519-618]). \\
\ni [Wie2] N. Wiener, {\sl The Fourier integral and certain of its applications}.  Cambridge Univ. Press, 1933. \\
\ni [Wie3] N. Wiener, {\sl Collected works with commentaries, Vol. 2.  Generalized harmonic analysis and Tauberian theory; classical harmonic and complex analysis}.  MIT Press, Cambridge MA, 1979. \\

\ni Mathematics Department, Imperial College, London SW7 2AZ; n.bingham@ic.ac.uk

\end{document}